\documentclass{article}

\usepackage{graphicx}
\usepackage{hyperref}

\begin{document}
\title{The Secrets of Notakto:  Winning at X-only Tic-Tac-Toe}

\author{Thane E. Plambeck, Greg Whitehead} 
\maketitle
\begin{abstract}
We analyze misere play of ``impartial'' tic-tac-toe---a game suggested by Bob Koca in which {\bf both}
players make X's on the board, and the first player to complete 
three-in-a-row {\bf loses}.  This game was recently discussed on 
\href{http://mathoverflow.net/questions/24693/neutral-tic-tac-toe}{mathoverflow.net} in a thread created by Timothy Y. Chow.  
\end{abstract}

\section{Introduction}

Suppose tic-toe-toe is played on the usual 3x3 board, but where both players make X's on the board. 
The first player to complete a line of three-in-a-row {\em loses} the game.  

Who should win?  The answer for a single 3x3 board is given in a recent 
\href{http://mathoverflow.net/questions/24693/neutral-tic-tac-toe}{mathoverflow.net discussion} \cite{chow}:

\begin{quotation}In the 3x3 misere game, the first player wins by playing in the center, and then wherever 
the second player plays, the first player plays a knight's move away from that.
\end{quotation}

\noindent Kevin Buzzard pointed out that any other first-player move loses:

\begin{quotation}
The reason any move other than the centre loses for [the first player to move] in the 3x3 game is that [the second player]
can respond with a move diametrically opposite [the first player's] initial move.   This makes the 
centre square unplayable, and then player two just plays the ``180 degree rotation'' strategy which clearly wins.
\end{quotation}

In this note we generalize these results to give a complete analysis of multiboard impartial tic-tac-toe under the 
disjunctive misere-play convention.

\section{Disjunctive misere play}

A {\em disjunctive} game of 3x3 impartial tic-tac-toe is played not just with one tic-tac-toe board, but more 
generally with an arbitrary (finite) number of such boards forming the start position. 
On a player's move, he or she selects a single one of the boards, and makes an X on it 
(a board that already has a three-in-a-row configuration of X's is considered unavailable for further moves and out of play). 

Play ends when every board has a three-in-a-row configuration. The player 
who completes the last three-in-a-row on the last available board is the loser.

\section{The misere quotient of 3x3 impartial tic-tac-toe}

We can give a succinct and complete analysis of the best misere play of an arbitrarily complicated disjunctive sum of impartial
3x3 tic-tac-toe positions by introducing a certain 18-element commutative monoid $Q$ given by the presentation 

\begin{equation}
\label{monoidrelations}
Q  =  \langle\ a,b,c,d\ | \ a^2=1,\ b^3=b,\ b^2c=c,\ c^3=ac^2,\ b^2d=d,\ cd=ad,\ d^2=c^2 \rangle .
\end{equation}

\noindent The monoid $Q$ has eighteen elements

\begin{equation}
\label{monoidelements}
Q = \{1, a, b, a b, b^2, a b^2, c, a c, b c, a b c, c^2, a c^2, b c^2, 
 a b c^2, d, a d, b d, a b d\},
\end{equation}

\noindent and it is called the {\bf misere quotient} of impartial tic-tac-toe\footnote{For the cognoscenti: 
$Q$ arises as the misere quotient of the hereditary closure of the sum $G$ of two impartial misere games $G = 4 + \{2+,0\}$.  
The game $\{2+,0\}$ is the misere canonical form of the 3x3 single board start position, 
and ``4'' represents the nim-heap of size 4, which also happens to occur as a single-board position in impartial tic-tac-toe.  
In describing these misere canonical forms, we've used the notation of John Conway's On Numbers and Games, 
on page 141, Figure 32.}.

A complete discussion of the misere quotient theory (and how $Q$ can be calculated from the rules of impartial tic-tac-toe) 
is outside the scope of this document.  General information about misere quotients and their construction
can be found in \cite{mq1}, \cite{mq2}, \cite{mq3}, and \cite{mq4}. 
One way to think of $Q$ is that it captures the misere analogue of the ``nimbers'' and ``nim addition'' that 
are used in normal play disjunctive impartial game analyses, but localized to the play of this 
particular impartial game, misere impartial 3x3 tic-tac-toe.

In the remainder of this paper, we simply take $Q$ as given.

\section{Outcome determination}

Figure~\ref{dictionary} (on page \pageref{dictionary}, after the References) assigns an element of $Q$ 
to each of the conceivable 102 non-isomorphic positions\footnote{We mean ``non-isomorphic'' under a reflection or rotation of the board.  
In making this count, we're including positions that couldn't be reached in actual play 
because they have too many completed rows of X's, but that doesn't matter since all 
those elements are assigned the identity element of Q.} in 3x3 single-board impartial 
tic-tac-toe. 

To determine the outcome of a multi-board position (ie, whether the position is an 
{\em N-position}---a Next player to move wins in best play, or alternatively, a {\em P-position}---second player to move wins), 
one first multiplies the corresponding elements of Q from the 
dictionary together.  The resulting word is then reduced via the relations \ref{monoidrelations}, that we 
started with above, necessarily eventually arriving at at one of the eighteen words (in the alphabet $a, b, c ,d$) that make
up the elements of $Q$.

If that word ends up being one of the four words in the set P

\begin{equation}
\label{ppositions}
P = \{a, b^2, bc, c^2\},
\end{equation}

\noindent the position is P-position; otherwise, it's an N-position.

\section{Example analysis}

\begin{figure}[h*]
\centering
\includegraphics[scale=0.63]{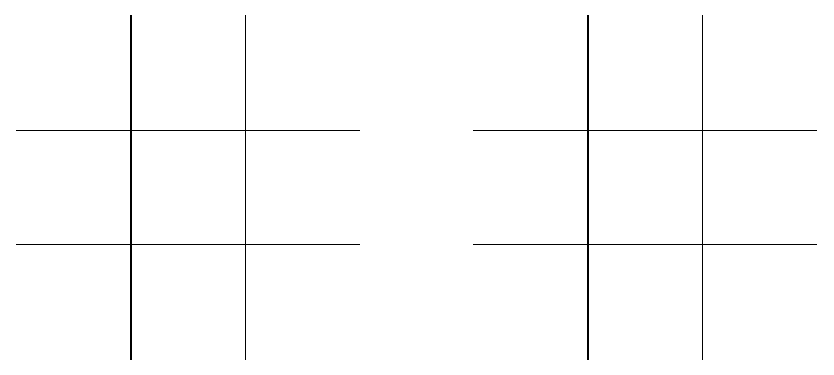}
\caption{\label{twoboardstart} The two-board start position.}
\end{figure}

To illustrate outcome calculation for Impartial Tic-Tac-Toe, we consider the two-board start position shown in Figure~\ref{twoboardstart}.

Consulting Figure~\ref{dictionary}, we find that the monoid-value of a single empty board is $c$.  Since we have two such boards in our position, we multiply these two values together
and obtain the monoid element
\[ c^2 = c \cdot c.\]

\begin{figure}[h*]
\centering
\includegraphics[scale=0.63]{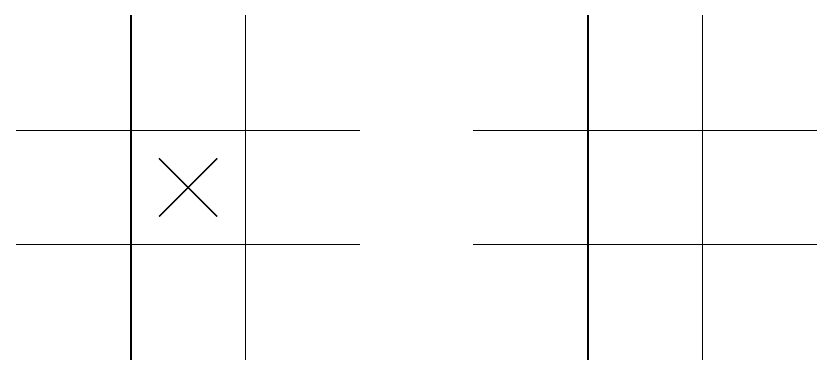}
\caption{\label{twoboardstartfirstmove} A doomed first move from the two-board start position.}
\end{figure}

Since $c^2$ is in the set P (equation (\ref{ppositions})), the position shown in Figure~\ref{twoboardstart} is a {\em second} player win.  Supposing therefore that we 
helpfully encourage our opponent to make the first
move, and that she moves to the center of one of the boards, we arrive at the position shown in Figure \ref{twoboardstartfirstmove}.

\begin{figure}[h*]
\centering
\includegraphics[scale=0.63]{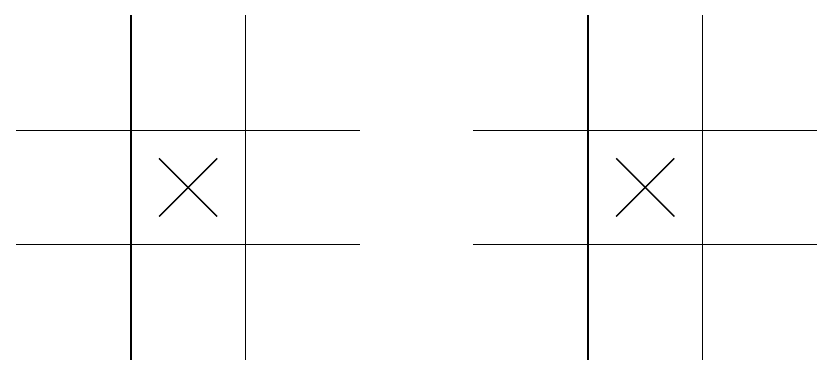}
\caption{\label{twoboardstartsecondmove} Mimicry works here, but not in general.}
\end{figure}

It so happens that if we mimic our opponent's move on the other board, this happens to be a winning move.
We arrive at the position shown in Figure \ref{twoboardstartsecondmove}, each of whose two boards is of value $c^2$;  
multiplying these two together, and simplifying via the relations shown in equation (\ref{monoidrelations}), we have

\begin{eqnarray*}
c^4 & = & c^3 \cdot c \\
    & = & a c^2 \cdot c \\
    & = & a c^3 \\
    & = & a a c^2 \\
    & = & c^2,
\end{eqnarray*}
which is a P-position, as desired.  

So is the general winning stategy of the two-board position simply to copy our opponent's moves on the other board?  Far from it:  consider what happens if our opponent should
decide to complete a line on one of the boards---copying that move on the other board, we'd {\bf lose} rather than win!  For example, from the N-position shown in Figure~\ref{nomimic},
there certainly is a winning move, but it's {\em not} to the upper-right-hand corner of the board on the right, which loses.

\begin{figure}[h*]
\centering
\includegraphics[scale=0.63]{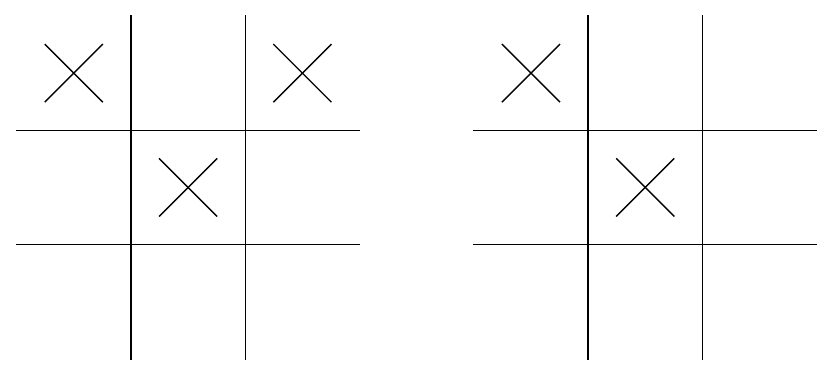}
\caption{\label{nomimic} An N-position in which mimickry loses.}
\end{figure}

We invite our reader to find a correct reply!

\section{The iPad game Notakto}

Evidently the computation of general outcomes in misere tic-tac-toe is somewhat complicated, involving 
computations in finite monoid and looking up values from a table of all possible single-board positions.

However, we've found that a human can develop the ability to 
win from multi-board positions with some practice.

\begin{figure}[h*]
\centering
\includegraphics[scale=0.63]{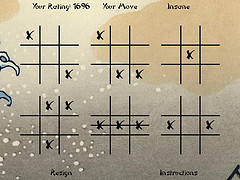}
\caption{\label{nomimic} A six-board game of Notakto, in progress.}
\end{figure}

{\bf Notakto} ``No tac toe'' is an iPad game that allows the user to practice playing misere X-only Tic-Tac-Toe against a computer.
Impartial misere tic-tac-toe from start positions involving one up to as many as six initial tic-tac-toe boards are supported.

The Notakto iPad application is available for free at {\tt http://www.notakto.com}.

\section{Final question}
Does the 4x4 game have a finite misere quotient? 

\begin{figure}[h*]
\centering
\includegraphics[scale=0.63]{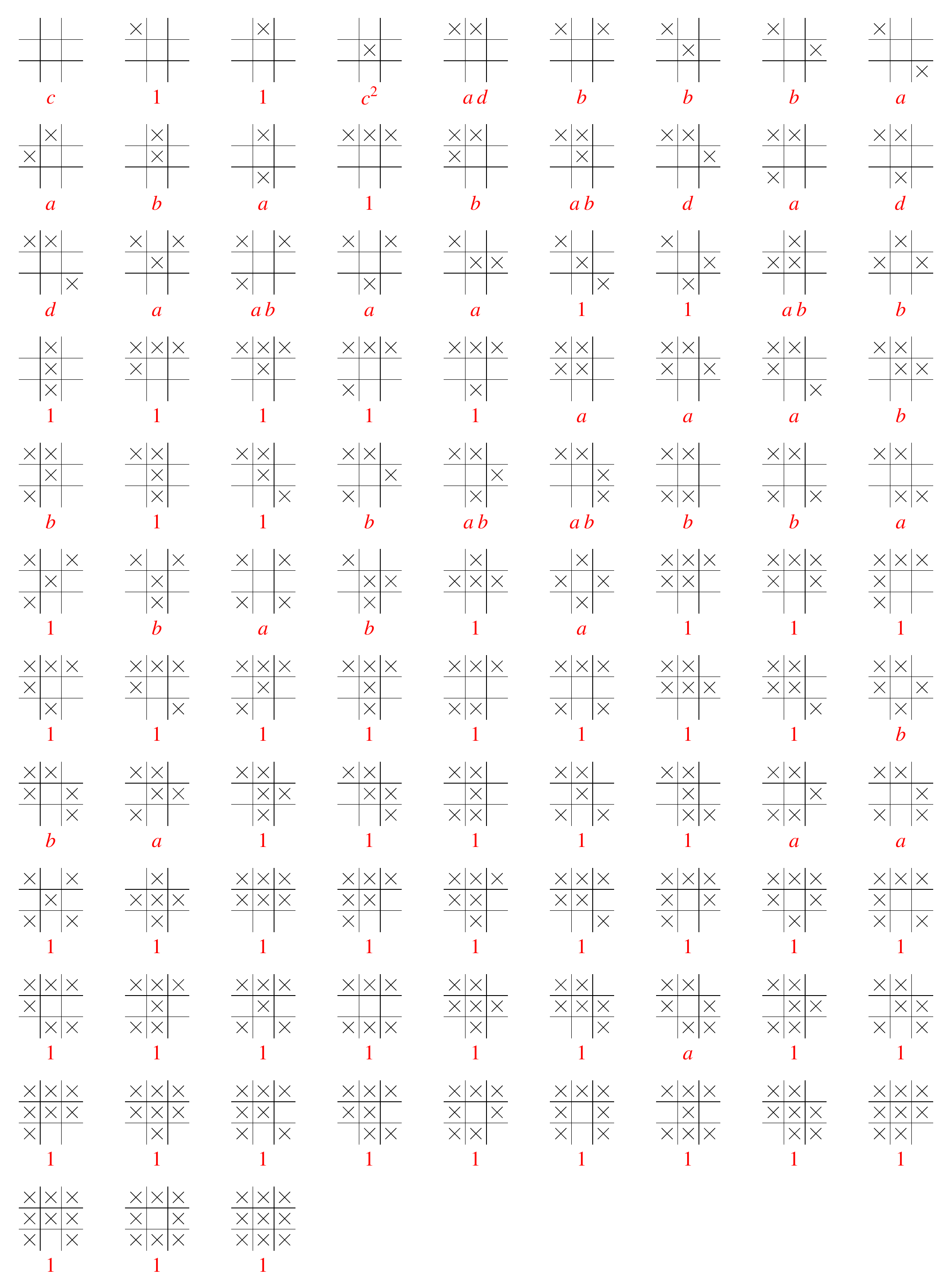}
\caption{\label{dictionary} The 102 nonisomorphic ways of arranging zero to nine X's on a tic-tac-toe board, each shown together with its corresponding misere quotient element from $Q$.}
\end{figure}

\end{document}